# Estimation in a class of nonlinear heteroscedastic time series models

### Joseph Ngatchou-Wandji

*Laboratoire de Mathématiques Nicolas Oresme,*
*CNRS-UMR 6139, Université de Caen, Campus II, BP 5186*
*Boulevard du Maréchal Juin, F-14032 Caen, France.*
*e-mail:* `joseph.ngatchou-wandji@unicaen.fr`

**Abstract:** Parameter estimation in a class of heteroscedastic time series models is investigated. The existence of conditional least-squares and conditional likelihood estimators is proved. Their consistency and their asymptotic normality are established. Kernel estimators of the noise's density and its derivatives are defined and shown to be uniformly consistent. A simulation experiment conducted shows that the estimators perform well for large sample size.

## Contents

## 1. Introduction

Although parametric models are prone to misspecification, they still be attractive because they describe concisely the link between the past observations and the predicted variable. Various parametric models have been proposed these last decades. For a review, see Brockwell and Davis (1991), Brockwell and Davis (1996), Shumway and Stoffer (2001), and Tong (1990). Parameter estimation for linear models has been widely studied, while for nonlinear models, because of their complexity, the study is done in general for tractable cases. There is an increasing interest in estimating the parameters of ARCH and GARCH models introduced respectively by Engle (1982) and Bollerslev (1986). Most of the existing literature assume a Gaussian error distribution and study the consistency and asymptotic normality of the conditional Gaussian likelihood estimators. Some relevant papers are Engle (1982), Weiss (1986) for ARCH models, and for ARCH($\infty$) and/or GARCH models, Bollerslev (1986), Lumsdaine (1996), Hyndman and Yao (2002), Francq and Zakoïan (2004), Robinson and Zaffroni (2006), Straumann and Mikosch (2006), and Francq and Zakoïan (2007). Other papers dealing with parameter estimation in heterocedastic models include the works of Giraitis and Robinson (2001) who propose a Wittle estimation for a class of parametric ARCH($\infty$), Chatterjee and Das (2003) who study estimators obtained by minimizing certain functionals for ARCH models, Peng and Yao (2003) who propose least absolute deviations estimators for ARCH and GARCH models, and Berkes and Horváth (2004) who study likelihood estimators for GARCH models.

In the present paper, we study parameter estimation for more general heterocedastic models. Precisely, we consider the class of identifiable parametric stochastic models

$$X_i = m\left(\rho; Z_{i-1}\right) + \sigma(\theta; Z_{i-1})\varepsilon_i, \ i \in \mathbb{Z}, \tag{1.1}$$

where $(X_i)_{i \in \mathbb{Z}}$ is stationary and ergodic; $(Z_i = (X_i, \ldots, X_{i-q+1}, X_{i-q}))_{i \in \mathbb{Z}}$ is a sequence of $q$-dimensional vector with $q$ being a nonnegative possibly infinite integer; $(\varepsilon_i)_{i \in \mathbb{Z}}$ is a sequence of iid centered random variables with unit variance such that $\varepsilon_i$ is independent of $\sigma(Z_j, j < i)$; the parameter column vector $\psi = (\rho', \theta')'$ belongs to $\Psi = \Theta \times \tilde{\Theta} \subset \mathbb{R}^I \times \mathbb{R}^J$, for some positive integers $I$ and $J$, and the functions $m\left(\rho; z\right)$ and $\sigma(\theta; z)$ have known forms. We aim to prove the existence of asymptotical normal estimators for the true parameter vector $\psi_0 = (\rho'_0, \theta'_0)'$, and uniformly consistent estimators for the noise's density and its derivatives, when this function exists.

The class of models (1.1) contains models such as ARMA, EXPAR, ARCH, GARCH, SETAR-ARCH, $\beta$-ARCH and many others. As far as the probabilist properties of these models are concerned, their invertibility is readily obtained for example for $|\sigma(\theta; z)| > 0$. For some of them (see, e.g., Ngatchou-Wandji (2005)), a sufficient condition for strict stationarity can be obtained e.g., by checking the conditions (S1)-(S4) of p. 86 in Taniguchi and Kakizawa (2000). The case of GARCH models which generalizes ARCH models has been studied by Chen and An (1998). Next, a sufficient condition for geometry ergodicity can be obtained by applying a result of Tjøstheim (1990), while for a particular class of ARCH models nested in (1.1), this property has been investigated by An, Chen and Huang (1997). Finally, it is possible that from the theory of Markov chains, other interesting conditions for stationarity and ergodicity be obtained for many models within (1.1).

Under mild conditions, a conditional least-squares estimator of $\rho_0$ is defined in McKeague and Zhang (1994). Its consistency and asymptotic normality is established. The same is done for $\theta_0$ in Ngatchou-Wandji (2002). Such results have also been established for multivariate nonlinear AR models by Tjøstheim (1986). Our main contribution is the study of the estimation of the couple of parameters $\psi_0 = (\rho_0', \theta_0')'$ in model (1.1) by conditional least-squares and conditional maximum likelihood methods, when the conditional distribution is non necessarily normal and $q$ possibly infinite. Our results generalize most of those based on least-squares and pseudo or quasi-likelihood estimation.

After the assumptions given in Section 2, we prove in Section 3 the existence of a sequence of asymptotical normal conditional least-squares estimators for $\psi_0$. Section 4 deals with the existence of conditional likelihood estimators for this parameter. In Section 5, we give some common examples comprised in (1.1). In Section 6, the estimation of the noise's density and its derivatives is investigated. A simulation study done in Section 7 ends our work.

## 2. General assumptions

In the whole text, the transpose of a vector or a matrix function $\mathcal{H}(x)$ is denoted by $\mathcal{H}'(x)$. Let $r$ be either $I$ or $J$. For given real functions $\mathcal{F}(\alpha; z)$ defined on a non-empty subset of $\mathbb{R}^r \times \mathbb{R}^q$ and $\mathcal{K}(\psi; z)$ defined on a non-empty subset of $\mathbb{R}^I \times \mathbb{R}^J \times \mathbb{R}^q$, we denote

$$\partial \mathcal{F}(\alpha; z) = \Big(\frac{\partial \mathcal{F}(\alpha; z)}{\partial \alpha_1}, \ldots, \frac{\partial \mathcal{F}(\alpha; z)}{\partial \alpha_r}\Big)', \quad \partial^2 \mathcal{F}(\alpha; z) = \Big(\frac{\partial^2 \mathcal{F}(\alpha; z)}{\partial \alpha_i \partial \alpha_j} : 1 \le i, j \le r\Big),$$

$$\partial_\rho \mathcal{K}(\psi; z) = \Big(\frac{\partial \mathcal{K}(\psi; z)}{\partial \rho_1}, \ldots, \frac{\partial \mathcal{K}(\psi; z)}{\partial \rho_I}\Big)', \quad \partial_\theta \mathcal{K}(\psi; z) = \Big(\frac{\partial \mathcal{K}(\psi; z)}{\partial \theta_1}, \ldots, \frac{\partial \mathcal{K}(\psi; z)}{\partial \theta_J}\Big)',$$

$$\partial^2_{\rho\theta} \mathcal{K}(\psi; z) = \Big(\frac{\partial^2 \mathcal{K}(\psi; z)}{\partial \rho_i \partial \theta_j} : 1 \le i \le I, 1 \le j \le J\Big),$$

$$\partial^2_{\theta\rho} \mathcal{K}(\psi; z) = \Big(\frac{\partial^2 \mathcal{K}(\psi; z)}{\partial \theta_j \partial \rho_i} : 1 \le j \le J, 1 \le i \le I\Big),$$

$$\partial^2_{\rho^2} \mathcal{K}(\psi; z) = \Big(\frac{\partial^2 \mathcal{K}(\psi; z)}{\partial \rho_i \partial \rho_j} : 1 \le i, j \le I\Big),$$

$$\partial^2_{\theta^2} \mathcal{K}(\psi; z) = \Big(\frac{\partial^2 \mathcal{K}(\psi; z)}{\partial \theta_i \partial \theta_j} : 1 \le i, j \le J\Big).$$

For a vector or matrix function $\mathcal{H}(x)$, we denote by $\partial' \mathcal{H}(x)$ the transpose of $\partial \mathcal{H}(x)$. With this, we define $\partial \mathcal{K}(\psi; z) = (\partial'_\rho \mathcal{K}(\psi; z); \partial'_\theta \mathcal{K}(\psi; z))'$. We also define

$$\partial^2 \mathcal{K}(\psi; z) = \begin{pmatrix} \partial^2_{\rho^2} \mathcal{K}(\psi; z) & \partial^2_{\rho\theta} \mathcal{K}(\psi; z) \\ \partial^2_{\theta\rho} \mathcal{K}(\psi; z) & \partial^2_{\theta^2} \mathcal{K}(\psi; z) \end{pmatrix}.$$

For a real-valued function $h$, $h^{(p)}$ denotes its $p$th order derivative, with $h^{(0)} = h$. We denote by $||\mathbf{V}||_{\mathcal{E}}$ the Euclidean norm of the vector $\mathbf{V}$ and by $||\mathbf{M}||_{\mathcal{M}} = \max_{i,j} |M_{ij}|$ the norm of the square matrix $\mathbf{M} = (M_{ij})$.

We next assume that the true parameter vector $\psi_0 = (\rho_0', \theta_0')'$ of (1.1) is such that $\rho_0 \in \text{int}(\Theta)$ and $\theta_0 \in \text{int}(\tilde{\Theta})$, where $\text{int}(\Theta)$ and $\text{int}(\tilde{\Theta})$ denote respectively the nonempty interior of $\Theta$ and $\tilde{\Theta}$. We also suppose that all the random variables in this paper are defined on the same probability space $(\Omega, \mathcal{W}, P)$, where $\Omega$ is a set, $\mathcal{W}$ a $\sigma$-field of $\Omega$ and $P$ a probability measure on $\mathcal{W}$. The following assumptions are needed:

$(\mathcal{A}_1)$ The common fourth order moment of the $\varepsilon_i$'s is finite.

$(\mathcal{A}_2)$ The functions $m(\rho; z)$ and $\sigma(\theta; z)$ are each twice continuously differentiable with respect to $\rho \in \mathrm{int}(\Theta)$ and $\theta \in \mathrm{int}(\tilde{\Theta})$ respectively, and there exists a positive function $\alpha(z)$ such that $E[\alpha^4(Z_0)] < \infty$ and

$$\max\{\sup_{\rho \in \Theta} |m(\rho; z)|, \sup_{\rho \in \Theta} ||\partial m(\rho; z)||_{\mathcal{E}}, \sup_{\rho \in \Theta} ||\partial^2 m(\rho; z)||_{\mathcal{M}},$$

$$\sup_{\theta \in \tilde{\Theta}} |\sigma(\theta; z)|, \sup_{\theta \in \tilde{\Theta}} ||\partial \sigma(\theta; z)||_{\mathcal{E}}, \sup_{\theta \in \tilde{\Theta}} ||\partial^2 \sigma(\theta; z)||_{\mathcal{M}}\} \le \alpha(z),$$

$(\mathcal{A}_3)$ There exists a positive function $\beta(z)$ such that $E[\beta^4(Z_0)] < \infty$ and for all $\rho_1, \rho_2 \in \Theta$ and $\theta_1, \theta_2 \in \tilde{\Theta}$,

$$\max\{|m(\rho_1; z) - m(\rho_2; z)|, ||\partial m(\rho_1; z) - \partial m(\rho_2; z)||_{\mathcal{E}},$$

$$||\partial^2 m(\rho_1; z) - \partial^2 m(\rho_2; z)||_{\mathcal{M}}, |\sigma(\theta_1; z) - \sigma(\theta_2; z)|,$$

$$||\partial \sigma(\theta_1; z) - \partial \sigma(\theta_2; z)||_{\mathcal{E}}, ||\partial^2 \sigma(\theta_1; z) - \partial^2 \sigma(\theta_2; z)||_{\mathcal{M}}\}$$

$$\le \beta(z) \min\{||\rho_1 - \rho_2||_{\mathcal{E}}, ||\theta_1 - \theta_2||_{\mathcal{E}}\}.$$

Assumption $(\mathcal{A}_1)$ is at least satisfied by Gaussian and Student $\varepsilon_i$'s. One can find in the literature, numbers of models with the functions $m(\rho; z)$ and $\sigma(\theta; z)$ satisfying $(\mathcal{A}_2)$ and $(\mathcal{A}_3)$ (see, e.g., Ngatchou-Wandji (2005)).

## 3. Conditional least-squares estimation

The purpose of this section is the study of the existence of estimators for $\psi_0 = (\rho_0', \theta_0')'$ by a conditional least-squares method. Recall that the conditional mean and the conditional variance functions of (1.1) are almost surely defined for all $z \in \mathbb{R}^q$ by $E(X_1 \mid Z_0 = z) = m(\rho; z)$ and $E\{[X_1 - m(\rho; Z_0)]^2 \mid Z_0 = z\} = \sigma^2(\theta; z)$. From these equalities, for any bounded measurable functions $\gamma(z)$ and $\lambda(z)$, we have $E[(X_1 - m(\rho; Z_0))\lambda(Z_0)] = 0$ and $E[\{(X_1 - m(\rho; Z_0))^2 - \sigma^2(\theta; Z_0)\}\gamma(Z_0)] = 0$. For estimating $\psi_0$, our idea is to search for the zeros of the gradients of the sample variances of the sequences of centered random variables $(X_i - m(\rho; Z_{i-1}))$, $i = 1, \ldots, n$ and $([X_i - m(\rho; Z_{i-1})]^2 - \sigma^2(\theta; Z_{i-1}))$, $i = 1, \ldots, n$.

Given $X_{-q}, \ldots, X_{-1}, X_0, X_1, \ldots, X_n$, denote $\mathbf{X}_n = (X_n, \ldots, X_1, X_0, X_{-1}, \ldots, X_{-q})$ and define the sequences of random functions

$$U_n(\rho; \mathbf{X}_n) = \sum_{i=1}^n [X_i - m(\rho; Z_{i-1})]^2 \lambda^2(Z_{i-1}) \tag{3.1}$$

$$S_n(\psi; \mathbf{X}_n) = \sum_{i=1}^n \left\{[X_i - m(\rho; Z_{i-1})]^2 - \sigma^2(\theta; Z_{i-1})\right\}^2 \gamma^2(Z_{i-1}) \tag{3.2}$$

and the matrices

$$\Phi_{11} = 2E[\lambda^2(Z_0)\sigma(\theta_0; Z_0)\partial m(\rho_0; Z_0)\partial' m(\rho_0; Z_0)]$$
$$\Phi_{22} = 8E[\gamma^2(Z_0)\sigma^2(\theta_0; Z_0)\partial\sigma(\theta_0; Z_0)\partial'\sigma(\theta_0; Z_0)],$$

assumed to be positive definite. Define also the matrices

$$\Delta_{11} = 4E[\lambda^4(Z_0)\sigma^2(\theta_0; Z_0)(\Phi_{11}^{-1})'\partial m(\rho_0; Z_0)\partial' m(\rho_0; Z_0)\Phi_{11}^{-1}]$$

$$\Delta_{12} = \Delta_{21}'$$
$$= 8E[\lambda^2(Z_0)\gamma^2(Z_0)\sigma^4(\theta_0; Z_0)(\Phi_{11}^{-1})'\partial m(\rho_0; Z_0)\partial'\sigma(\theta_0; Z_0)\Phi_{22}^{-1}]E[\varepsilon_0(\varepsilon_0^2 - 1)]$$

$$\Delta_{22} = 16E[\gamma^4(Z_0)\sigma^6(\theta_0; Z_0)(\Phi_{22}^{-1})'\partial\sigma(\theta_0; Z_0)\partial'\sigma(\theta_0; Z_0)\Phi_{22}^{-1}]E[(\varepsilon_0^2 - 1)^2],$$

and

$$\Delta = \left( \begin{array}{cc} \Delta_{11} & \Delta_{12} \\ \Delta_{21} & \Delta_{22} \end{array} \right).$$

**Theorem 3.1.** *Assume that the assumptions $(\mathcal{A}_1)$-$(\mathcal{A}_3)$ hold and $\Delta$ is positive definite. Then,*

(i) *there exists a sequence of estimators $\tilde{\psi}_n = (\tilde{\rho}_n', \tilde{\theta}_n')'$ such that $\tilde{\psi}_n \xrightarrow{a.s.} \psi_0$, and for any $\epsilon > 0$, there exists an event $\mathcal{S}_1$ with $P(\mathcal{S}_1) > 1 - \epsilon$, and a nonnegative integer $n_1$ such that on $\mathcal{S}_1$, for $n > n_1$,*

   - $\partial U_n(\tilde{\rho}_n; \mathbf{X}_n) = \mathbf{0}$ *and $U_n(\rho; \mathbf{X}_n)$ attains a relative minimum at $\rho = \tilde{\rho}_n$*
   - • *assuming $\tilde{\rho}_n$ fixed, $\partial_\theta S_n((\tilde{\rho}_n, \tilde{\theta}_n); \mathbf{X}_n) = \mathbf{0}$ and $S_n((\tilde{\rho}_n, \theta); \mathbf{X}_n)$ attains a relative minimum at $\theta = \tilde{\theta}_n$.*

(ii) $n^{1/2}(\tilde{\psi}_n - \psi_0) \xrightarrow{\mathcal{D}} \mathcal{N}(\mathbf{0}, \Delta).$

*Proof.* It suffices to check the hypotheses of Theorem 3.2.23 of Taniguchi and Kakizawa (2000), established by Klimko and Nelson (1978) by using Egorov Theorem (see, e.g., Taniguchi and Kakizawa (2000), p. 97). From simple computations one obtains:

$$\partial U_n(\rho; \mathbf{X}_n) = -2\sum_{i=1}^{n} \lambda^2(Z_{i-1})\partial m(\rho; Z_{i-1})\left[X_i - m(\rho; Z_{i-1})\right]$$

and

$$\partial^2 U_n(\rho; \mathbf{X}_n)$$
$$= 2\sum_{i=1}^{n} \lambda^2(Z_{i-1})\Big(\partial m(\rho; Z_{i-1})\partial' m(\rho; Z_{i-1}) - \partial^2 m(\rho; Z_{i-1})\left[X_i - m(\rho; Z_{i-1})\right]\Big).$$

By ergodicity, it is immediate that, as $n$ tends to infinity,

$$\frac{1}{n}\partial U_n(\rho_0; \mathbf{X}_n) \xrightarrow{a.s.} \mathbf{0} \quad \text{and} \quad \frac{1}{n}\partial^2 U_n(\rho_0; \mathbf{X}_n) \xrightarrow{a.s.} \Phi_{11}.$$

For any vector $\rho_* \in \text{int}(\Theta)$ define the sequence of random matrix functions

$$V_n(\rho_*; \mathbf{X}_n) = \partial^2 U_n(\rho_*; \mathbf{X}_n) - \partial^2 U_n(\rho_0; \mathbf{X}_n),$$

and denote by $V_n(\rho_*; \mathbf{X}_n)_{\ell k}$ its $(\ell, k)th$ entry. Then,

$$V_n(\rho_*; \mathbf{X}_n)_{\ell k}$$
$$= 2 \sum_{i=1}^{n} \lambda^2(Z_{i-1}) \left\{ \left( \frac{\partial m(\rho_*; Z_{i-1})}{\partial \rho_\ell} \frac{\partial m(\rho_*; Z_{i-1})}{\partial \rho_k} - \frac{\partial m(\rho_0; Z_{i-1})}{\partial \rho_\ell} \frac{\partial m(\rho_0; Z_{i-1})}{\partial \rho_k} \right) \right.$$
$$\left. - \left( \partial^2 m(\rho_*; Z_{i-1}) \left[ X_i - m(\rho_*; Z_{i-1}) \right] - \partial^2 m(\rho_0; Z_{i-1}) \left[ X_i - m(\rho_0; Z_{i-1}) \right] \right) \right\}.$$

Thus, in view of $(\mathcal{A}_1)$-$(\mathcal{A}_3)$, it is easy to see that there exists a positive real-valued function $\nu_{\ell k}(z)$ with $E[\nu_{\ell k}^4(Z_0)] < \infty$ such that

$$|V_n(\rho_*; \mathbf{X}_n)_{\ell k}| \leq ||\rho_* - \rho_0||_{\mathcal{E}} \sum_{i=1}^{n} \nu_{\ell k}(Z_i).$$

Now for $\delta > 0$ such that $||\rho - \rho_0||_{\mathcal{E}} < \delta$, and for $\rho_*$ lying between $\rho$ and $\rho_0$, we have by the above inequality that:

$$\frac{1}{n\delta} |V_n(\rho_*; \mathbf{X}_n)_{\ell k}| \leq \frac{1}{n\delta} ||\rho_* - \rho_0||_{\mathcal{E}} \sum_{i=1}^{n} \nu_{\ell k}(Z_{i-1}) \leq \frac{1}{n} \sum_{i=1}^{n} \nu_{\ell k}(Z_{i-1}).$$

Next, by ergodicity, the right-hand side of the last inequality converges a.s. to $E[\nu_{\ell k}(Z_0)] < \infty$ as $n$ tends to infinity. It is then clear that for any $(\ell, k)$,

$$\lim_{n \to \infty} \sup_{\delta \to 0} \frac{1}{n\delta} |V_n(\rho_*; \mathbf{X}_n)_{\ell k}| < \infty. \tag{3.3}$$

From Theorem 3.2.23 of Taniguchi and Kakizawa (2000), it follows that there exists a sequence of estimators $\tilde{\rho}_n$ such that $\tilde{\rho}_n \longrightarrow \rho_0$ almost surely, as $n \to \infty$ and for $\epsilon > 0$, one can find an event $E_1$ with $P(E_1) > 1 - \epsilon$ and a nonnegative integer $\tilde{n}$ such that on $E_1$, for $n > \tilde{n}$, $\partial U_n(\tilde{\rho}_n; \mathbf{X}_n) = \mathbf{0}$ and $U_n(\rho; \mathbf{X}_n)$ attains a relative minimum at $\rho = \tilde{\rho}_n$. The first part of $(i)$ is then handled. For the second part, for fixed $\tilde{\rho}_n$, we have from simple computations:

$$\partial_\theta S_n((\tilde{\rho}_n, \theta_0); \mathbf{X}_n)$$
$$= -4 \sum_{i=1}^{n} \gamma^2(Z_{i-1}) \sigma(\theta_0; Z_{i-1}) \partial \sigma(\theta_0; Z_{i-1})$$
$$\times \left\{ [X_i - m(\tilde{\rho}_n; Z_{i-1})]^2 - \sigma^2(\theta_0; Z_i) \right\}$$
$$= -4 \sum_{i=1}^{n} \gamma^2(Z_{i-1}) \sigma(\theta_0; Z_{i-1}) \partial \sigma(\theta_0; Z_{i-1})$$
$$\times \left\{ \sigma^2(\theta_0; Z_{i-1})(\varepsilon_i^2 - 1) + 2\sigma(\theta_0; Z_{i-1}) \varepsilon_i [m(\rho_0; Z_{i-1}) - m(\tilde{\rho}_n; Z_{i-1})] \right.$$
$$\left. + [m(\rho_0; Z_i) - m(\tilde{\rho}_n; Z_{i-1})]^2 \right\} \tag{3.4}$$

and

$$\partial_{\theta^2}^2 S_n((\tilde{\rho}_n, \theta_0); \mathbf{X}_n)$$

$$= 8 \sum_{i=1}^n \gamma^2(Z_{i-1}) \sigma^2(\theta_0; Z_{i-1}) \partial \sigma(\theta_0; Z_{i-1}) \partial' \sigma(\theta_0; Z_{i-1})$$

$$- 4 \sum_{i=1}^n \gamma^2(Z_{i-1}) \Big\{ \Big( \partial \sigma(\theta_0; Z_{i-1}) \partial' \sigma(\theta_0; Z_{i-1}) + \sigma(\theta_0; Z_{i-1}) \partial^2 \sigma(\theta_0; Z_{i-1}) \Big)$$

$$\times \Big( [X_i - m(\tilde{\rho}_n; Z_{i-1})]^2 - \sigma^2(\theta_0; Z_{i-1}) \Big) \Big\}$$

$$= 8 \sum_{i=1}^n \gamma^2(Z_{i-1}) \sigma(\theta_0; Z_{i-1}) \partial \sigma(\theta_0; Z_{i-1}) \partial' \sigma(\theta_0; Z_{i-1})$$

$$- 4 \sum_{i=1}^n \gamma^2(Z_{i-1}) \Big( \partial \sigma(\theta_0; Z_{i-1}) \partial' \sigma(\theta_0; Z_{i-1}) + \sigma(\theta_0; Z_{i-1}) \partial^2 \sigma(\theta_0; Z_{i-1}) \Big)$$

$$\times \Big\{ \sigma^2(\theta_0; Z_{i-1})(\varepsilon_i^2 - 1) + 2\sigma(\theta_0; Z_{i-1})\varepsilon_i [m(\rho_0; Z_{i-1}) - m(\tilde{\rho}_n; Z_{i-1})]$$

$$+ [m(\rho_0; Z_{i-1}) - m(\tilde{\rho}_n; Z_{i-1})]^2 \Big\}. \tag{3.5}$$

In view of $(\mathcal{A}_1)$-$(\mathcal{A}_3)$, applying the mean value theorem to (3.4) and (3.5), it is clear by ergodicity that as $n$ tends to infinity,

$$\frac{1}{n} \partial_\theta S_n((\tilde{\rho}_n, \theta_0); \mathbf{X}_n) \xrightarrow{a.s.} \mathbf{0} \quad \text{and} \quad \frac{1}{n} \partial_{\theta^2}^2 S_n((\tilde{\rho}_n, \theta_0); \mathbf{X}_n) \xrightarrow{a.s.} \Phi_{22}.$$

For any vector $\theta_* \in \text{int}(\tilde{\Theta})$ define the sequence of random functions

$$T_n(\theta_*; \mathbf{X}_n) = \partial_{\theta^2}^2 S_n(\theta_*; \mathbf{X}_n) - \partial_{\theta^2}^2 S_n(\theta_0; \mathbf{X}_n),$$

and denote by $T_n(\theta_*; \mathbf{X}_n)_{\ell k}$ its $(\ell, k)th$ entry.

$$T_n(\theta_*; \mathbf{X}_n)_{\ell k}$$

$$= 8 \sum_{i=1}^n \gamma^2(Z_{i-1}) \Big\{ \sigma(\theta_*; Z_{i-1}) \frac{\partial \sigma(\theta_*; Z_{i-1})}{\partial \theta_\ell} \frac{\partial \sigma(\theta_*; Z_{i-1})}{\partial \theta_k}$$

$$- \sigma(\theta_0; Z_{i-1}) \frac{\partial \sigma(\theta_0; Z_{i-1})}{\partial \theta_\ell} \frac{\partial \sigma(\theta_0; Z_{i-1})}{\partial \theta_k} \Big\}$$

$$- 4 \sum_{i=1}^n \gamma^2(Z_{i-1}) \Big\{ \Big[ \Big( \frac{\partial \sigma(\theta_*; Z_{i-1})}{\partial \theta_\ell} \frac{\partial \sigma(\theta_*; Z_{i-1})}{\partial \theta_k} + \sigma(\theta_*; Z_{i-1}) \frac{\partial^2 \sigma(\theta_*; Z_{i-1})}{\partial \theta_\ell \partial \theta_k} \Big)$$

$$\times \Big( [X_i - m(\tilde{\rho}_n; Z_{i-1})]^2 - \sigma^2(\theta_*; Z_{i-1}) \Big)$$

$$- \Big( \frac{\partial \sigma(\theta_0; Z_{i-1})}{\partial \theta_\ell} \frac{\partial \sigma(\theta_0; Z_{i-1})}{\partial \theta_k} + \sigma(\theta_0; Z_{i-1}) \frac{\partial^2 \sigma(\theta_0; Z_{i-1})}{\partial \theta_\ell \partial \theta_k} \Big)$$

$$\times \Big( [X_i - m(\tilde{\rho}_n; Z_{i-1})]^2 - \sigma^2(\theta_0; Z_{i-1}) \Big) \Big] \Big\}.$$

In view of $(\mathcal{A}_1)$-$(\mathcal{A}_3)$, it is easy to see that there exists a positive real-valued function $\varrho_{\ell k}(z)$ with $E[\varrho_{\ell k}^4(Z_0)] < \infty$ such that

$$|T_n(\theta_*; \mathbf{X}_n)_{\ell k}| \leq ||\theta_* - \theta_0||_{\mathcal{E}} \sum_{i=1}^n \varrho_{\ell k}(Z_{i-1}).$$

Again, for $\delta > 0$ such that $||\theta - \theta_0||_{\mathcal{E}} < \delta$, and for $\theta_*$ lying between $\theta$ and $\theta_0$, we have from above that:

$$\frac{1}{n\delta}|T_n(\theta_*; \mathbf{X}_n)_{\ell k}| \leq \frac{1}{n\delta}||\theta_* - \theta_0||_{\mathcal{E}} \sum_{i=1}^n \varrho_{\ell k}(Z_{i-1}) \leq \frac{1}{n}\sum_{i=1}^n \varrho_{\ell k}(Z_{i-1}).$$

It is easy to see that by the ergodic theorem, the right-hand side of the last inequality converges almost surely to $E[\varrho_{\ell k}(Z_0)] < \infty$ as $n$ tends to infinity. It is then clear that for any $(\ell, k)$, (3.3) holds with $T_n(\theta_*; \mathbf{X}_n)_{\ell k}$. Whence, applying Theorem 3.2.23 of Taniguchi and Kakizawa (2000), one can find a sequence of estimators $\tilde{\theta}_n$ such that $\tilde{\theta}_n \longrightarrow \theta_0$ almost surely, as $n \to \infty$ and for $\epsilon > 0$, one can find an event $E_2$ with $P(E_2) > 1 - \epsilon$ and a nonnegative integer $\hat{n}$ such that on $E_1 \cap E_2$, for $n > \hat{n}$, $\partial_\theta S_n(\tilde{\psi}_n; \mathbf{X}_n) = \mathbf{0}$ and $S_n((\tilde{\rho}_n, \theta); \mathbf{X}_n)$ attains a relative minimum at $\theta = \tilde{\theta}_n$. It is an easy matter to see that for all $\epsilon > 0$, $P(E_1 \cap E_2) > 1 - \epsilon$. Thus taking $\mathcal{S}_1 = E_1 \cap E_2$ and $n_1 = \max(\tilde{n}, \hat{n})$ yields the first part of Theorem 3.1. To handle the second point we observe that

$$\frac{1}{\sqrt{n}}\partial U_n(\rho_0; \mathbf{X}_n) = -\frac{2}{\sqrt{n}}\sum_{i=1}^n \lambda^2(Z_{i-1})\partial m(\rho_0; Z_{i-1})\sigma(\theta_0; Z_{i-1})\varepsilon_i,$$

and by a Taylor expansion of order one of the function $\partial U_n(\rho; \mathbf{X}_n)$ around $\rho_0$, for larger values of $n$, one can write

$$\sqrt{n}(\tilde{\rho}_n - \rho_0) = \frac{2}{\sqrt{n}}\sum_{i=1}^n \lambda^2(Z_{i-1})\sigma(\theta_0; Z_i)\varepsilon_i \partial' m(\rho_0; Z_{i-1})\Phi_{11}^{-1} + o_P(1).$$

One can also observe that

$$\frac{1}{\sqrt{n}}\partial_\theta S_n((\tilde{\rho}_n, \theta_0); \mathbf{X}_n) = -\frac{4}{\sqrt{n}}\sum_{i=1}^n \gamma^2(Z_{i-1})\sigma^3(\theta_0; Z_{i-1})\partial\sigma(\theta_0; Z_{i-1})(\varepsilon_i^2 - 1) + o_P(1)$$

and write for larger values of $n$,

$$\sqrt{n}(\tilde{\theta}_n - \theta_0) = \frac{4}{\sqrt{n}}\sum_{i=1}^n \gamma^2(Z_{i-1})\sigma^3(\theta_0; Z_{i-1})(\varepsilon_i^2 - 1)\partial'\sigma(\theta_0; Z_{i-1})\Phi_{22}^{-1} + o_P(1).$$

Then putting in Theorem 1 of Ngatchou-Wandji (2005): $\omega_i = \varepsilon_i$, $Y_i = Z_{i-1}$, $\Gamma_1(x) = x$, $\Gamma_2(x) = x^2 - 1$, $\Pi_1(z) = 2\lambda^2(z)\sigma(\theta_0; z) \partial' m(\rho_0; z)\Phi_{11}^{-1}$, and $\Pi_2(z) = 4\gamma^2(z)\sigma^3(\theta_0; z) \partial'\sigma(\theta_0; z)\Phi_{22}^{-1}$, it results that

$$\sqrt{n}(\tilde{\psi}_n - \psi_0) \xrightarrow{\mathcal{D}} \mathcal{N}(\mathbf{0}, \Delta).$$

$\square$

**Corollary 3.1.** *Assume that the assumptions of Theorem 3.1 hold and that $E[\varepsilon_0(\varepsilon_0^2 - 1)] = 0$. Then $\tilde{\rho}_n$ and $\tilde{\theta}_n$ are asymptotically uncorrelated.*

**Remark 3.1.** *The condition $E[\varepsilon_0(\varepsilon_0^2 - 1)] = E(\varepsilon_0^3) = 0$ in the above corollary holds for symmetric densities. If it does not hold, $\tilde{\rho}_n$ and $\tilde{\theta}_n$ will not be independent in general. This fact is ignored when the estimation is done for example with Gaussian or Student $\varepsilon_i$'s.*

**Remark 3.2.** *One could also prove the existence of consistent conditional type estimators for $\psi_0$ by using directly the random function $S_n(\psi; \mathbf{X}_n)$. It is clear that this would have provided a one step estimator. However, one may not retrieve the classical least-squares estimators. For example, in the very simple case of $m(\rho, z) = \rho z$ and $\sigma(\theta; z) = 1$, it is very difficult to have a simple expression for $\tilde{\rho}_n$ by minimizing $S_n((\rho, 1); \mathbf{X}_n)$, whereas, the preceding two-steps method yields the traditional least-squares estimator of $\rho_0$. Another approach (see, e.g., Heyde (1997)) consists in minimizing the sum of square $\sum_{i=1}^n \{[X_i - m(\rho; Z_{i-1})]/\sigma(\theta; Z_{i-1})\}^2$. Yet, for the special case of $\sigma(\theta; z) = \theta \neq 0$ and $\rho \in \mathbb{R}$, it is not clear how to estimate $\theta_0$ when $\rho_0 \neq \theta_0$.*

**Remark 3.3.** *The choice of the functions $\gamma(z)$ and $\lambda(z)$ is an open problem and we will not try to tackle it here. In the simulation, they are taken constant. Although this choice may be sub-optimal, it matches with what is done in the literature.*

## 4. Conditional likelihood estimation

For models such as (1.1), the density function of the noise can be useful for writing the likelihood and/or conditional likelihood functions. In practice, for choosing this density function, (1.1) can first be fitted by least-squares methods. Next, various tests can be applied to the residuals from the fitted model to help postulating an adequate density function $f$ (not necessarily Gaussian) for the noise. However, because it facilitates parameters estimation, pseudo-likelihood estimation method is very popular in practice. This probably explains the huge literature on the subject (see references given in Section 1).

In this section, we study the conditional likelihood estimation of the parameters when the noise has a non necessarily Gaussian density function $f$. This work has been done by Berkes and Horváth (2004) in the case of GARCH models. For simplicity, we restrict our study to models (1.1) for which the function $\sigma(\theta; z)$ satisfies:

$(\mathcal{B}_0)$ For all $(\theta, z) \in \mathbb{R}^J \times \mathbb{R}^q$, $\sigma(\theta; z) \geq \kappa$, for some constant $\kappa > 0$.

Under $(\mathcal{B}_0)$, the log-likelihood of a given sample $\mathbf{X}_n = (X_n, \ldots, X_1, X_0, X_{-1}, \ldots, X_{-q})$ conditional to $Z_0$ is

$$L_n(\psi; \mathbf{X}_n) = \sum_{i=1}^n \left\{ -\log[\sigma(\theta; Z_{i-1})] + \log[f(\varepsilon_i(\psi))] \right\}, \qquad (4.1)$$

where we recall that $\psi = (\rho', \theta')' \in \Psi$ and for all $i \in \mathbb{Z}$, $\varepsilon_i(\psi) = [X_i - m(\rho; Z_{i-1})]/\sigma(\theta; Z_{i-1})$. For the derivation of the results of this section, we make the following assumptions on the density function $f$:

$(\mathcal{B}_1)$ $f(x) > 0$ for all $x \in \mathbb{R}$, and $f$ is twice differentiable.

$(\mathcal{B}_2)$ $\phi_f = -f^{(1)}/f$ is differentiable with continuous derivative.

Next, for all $i \in \mathbb{Z}$, we define on $\mathbb{R}^I \times \mathbb{R}^J$, the following random functions:

$$
\begin{aligned}
\xi_i(\psi) &= \phi_f(\varepsilon_i(\psi)) \\
\dot{\xi}_i(\psi) &= \phi_f^{(1)}(\varepsilon_i(\psi)) \\
\zeta_i(\psi) &= \varepsilon_i(\psi)\phi_f(\varepsilon_i(\psi)) \\
\dot{\zeta}_i(\psi) &= \varepsilon_i(\psi)\phi_f^{(1)}(\varepsilon_i(\psi)) \\
\ddot{\zeta}_i(\psi) &= \zeta_i(\psi) + \varepsilon_i(\psi)\dot{\zeta}_i(\psi).
\end{aligned}
$$

We also need the following additional requirements:

$(\mathcal{B}_3)$ There exist a positive function $v(z)$ such that $E[v^4(Z_0)] < \infty$ and for all $i \in \mathbb{Z}$ and $\psi_1, \psi_2 \in \Psi$, a.s.,

$$
\begin{aligned}
\max\{&|\xi_i(\psi_1) - \xi_i(\psi_2)|, |\dot{\xi}_i(\psi_1) - \dot{\xi}_i(\psi_2)|, |\zeta_i(\psi_1) - \zeta_i(\psi_2)| \\
&|\dot{\zeta}_i(\psi_1) - \dot{\zeta}_i(\psi_2)|, |\ddot{\zeta}_i(\psi_1) - \ddot{\zeta}_i(\psi_2)|\} \le v(Z_i)\|\psi_1 - \psi_2\|_\varepsilon.
\end{aligned}
$$

$(\mathcal{B}_4)$ There exists a positive function $\tau(z)$, such that $E[\tau^4(Z_0)] < \infty$ and for all $i \in \mathbb{Z}$, a.s.

$$
\max\{\sup_{\psi \in \Psi}|\xi_i(\psi)|, \sup_{\psi \in \Psi}|\dot{\xi}_i(\psi)|, \sup_{\psi \in \Psi}|\zeta_i(\psi)|, \sup_{\psi \in \Psi}|\dot{\zeta}_i(\psi)|, \sup_{\psi \in \Psi}|\ddot{\zeta}_i(\psi)|\} \le \tau(Z_i).
$$

Such assumptions have been done in Ngatchou-Wandji (2005). They are at least satisfied by linear autoregressive models, EXPAR and TAR models, ARCH and more generally $\beta$-ARCH models, with Gaussian $f$.

Define the matrices

$$
\Sigma_{11} = E[\sigma^{-2}(\theta_0; Z_0)\partial m(\rho_0; Z_0)\partial' m(\rho_0; Z_0)] \int \phi_f^{(1)}(x)f(x)dx
$$

$$
\Sigma_{12} = \Sigma_{21}' = E[\sigma^{-2}(\theta_0; Z_0)\partial m(\rho_0; Z_0)\partial' \sigma(\theta_0; Z_0)] \int x\phi_f^2(x)f(x)dx
$$

$$
\Sigma_{22} = E[\sigma^{-2}(\theta_0; Z_0)\partial \sigma(\rho_0; Z_0)\partial' \sigma(\theta_0; Z_0)] \int x(\phi_f(x) + x\phi_f^{(1)}(x))f(x)dx
$$

$$
\Lambda_{11} = E[\sigma^{-2}(\theta_0; Z_0)\partial m(\rho_0; Z_0)\partial' m(\rho_0; Z_0)] \int \phi_f^2(x)f(x)dx
$$

$$
\Lambda_{12} = \Lambda_{21}' = E[\sigma^{-2}(\theta_0; Z_0)\partial \sigma(\theta_0; Z_0)\partial' m(\rho_0; Z_0)] \int \phi_f(x)(x\phi_f(x) - 1)f(x)dx
$$

$$
\Lambda_{22} = E[\sigma^{-2}(\theta_0; Z_0)\partial \sigma(\rho_0; Z_0)\partial' \sigma(\theta_0; Z_0)] \int (x\phi_f(x) - 1)^2 f(x)dx,
$$

$$\Sigma = \left( \begin{array}{cc} \Sigma_{11} & \Sigma_{12} \\ \Sigma_{21} & \Sigma_{22} \end{array} \right) \quad \text{and} \quad \Lambda = \left( \begin{array}{cc} \Lambda_{11} & \Lambda_{12} \\ \Lambda_{21} & \Lambda_{22} \end{array} \right).$$

**Theorem 4.1.** *Assume that* $(\mathcal{A}_1)$-$(\mathcal{A}_3)$ *and* $(\mathcal{B}_0)$-$(\mathcal{B}_4)$ *hold, and that the matrix* $\Sigma$ *is negative definite. If* $\int \phi_f(x) f(x) dx = 0$ *and* $\int x \phi_f(x) f(x) dx = 1$, *then*

(i) *there exists a sequence of estimators* $\hat{\psi}_n = (\hat{\rho}'_n, \hat{\theta}'_n)'$ *such that* $\hat{\psi}_n \xrightarrow{a.s.} \psi_0$, *and for any* $\epsilon > 0$, *there exists an event* $\mathcal{S}_2$ *with* $P(\mathcal{S}_2) > 1 - \epsilon$, *and a nonnegative integer* $n_2$ *such that on* $\mathcal{S}_2$, *for* $n > n_2$, $\partial L_n(\hat{\psi}_n; \mathbf{X}_n) = \mathbf{0}$ *and* $L_n(\psi; \mathbf{X}_n)$ *attains a relative maximum at* $\psi = \hat{\psi}_n$.

(ii) $n^{1/2}(\hat{\psi}_n - \psi_0) \xrightarrow{\mathcal{D}} \mathcal{N}(\mathbf{0}, \Sigma^{-1} \Lambda \Sigma^{-1})$.

*Proof.* The tools for the proof are exactly the same as those of the proof of Theorem 3.1. Define $Q_n(\psi; \mathbf{X}_n) = -L_n(\psi; \mathbf{X}_n)$. Then

$$\begin{aligned} &\partial Q_n(\psi_0; \mathbf{X}_n) \\ &= -\sum_{i=1}^{n} \sigma^{-1}(\theta_0; Z_{i-1}) \Big( \partial' m(\rho_0; Z_i) \xi_i(\psi_0), \ \partial' \sigma(\theta_0; Z_{i-1})(\zeta_i(\psi_0) - 1) \Big)', \end{aligned}$$

$$\begin{aligned} \partial^2_{\rho^2} Q_n(\psi; \mathbf{X}_n) &= -\sum_{i=1}^{n} \sigma^{-1}(\theta_0; Z_{i-1}) \Big( \partial^2 m(\rho_0; Z_{i-1}) \xi_i(\psi_0) \\ &\quad - \sigma^{-1}(\theta_0; Z_{i-1}) \partial m(\rho_0; Z_{i-1}) \partial' m(\rho_0; Z_{i-1}) \dot{\xi}_i(\psi_0) \Big), \end{aligned}$$

$$\partial^2_{\theta_0 \rho} Q_n(\psi_0; \mathbf{X}_n) = \sum_{i=1}^{n} \sigma^{-2}(\theta_0; Z_{i-1}) \{\xi_i(\psi_0) + \dot{\zeta}_i(\psi_0)\} \partial m(\rho_0; Z_{i-1}) \partial' \sigma(\theta_0; Z_{i-1})$$

$$\partial^2_{\rho\theta} Q_n(\psi_0; \mathbf{X}_n) = \sum_{i=1}^{n} \sigma^{-2}(\theta_0; Z_{i-1}) \{\xi_i(\psi_0) + \dot{\zeta}_i(\psi_0)\} \partial \sigma(\theta_0; Z_{i-1}) \partial' m(\rho_0; Z_{i-1}),$$

$$\begin{aligned} &\partial^2_{\theta^2} Q_n(\psi_0; \mathbf{X}_n) \\ &= \sum_{i=1}^{n} \Big[ \sigma^{-1}(\theta_0; Z_{i-1}) \Big\{ \sigma^{-1}(\theta_0; Z_{i-1}) \partial \sigma(\theta_0; Z_{i-1}) \partial' \sigma(\theta_0; Z_{i-1}) \\ &\quad - \partial^2 \sigma(\theta_0; Z_{i-1}) \Big\} (\zeta_i(\psi_0) - 1) + \sigma^{-2}(\theta_0; Z_{i-1}) \partial \sigma(\theta_0; Z_{i-1}) \partial' \sigma(\theta_0; Z_{i-1}) \ddot{\zeta}_i(\psi_0) \Big]. \end{aligned}$$

It is easy to see that, as $n$ tends to infinity,

$$\frac{1}{n} \partial Q_n(\psi_0; \mathbf{X}_n) \xrightarrow{a.s.} \mathbf{0} \quad \text{and} \quad \frac{1}{n} \partial^2 Q_n(\psi_0; \mathbf{X}_n) \xrightarrow{a.s.} -\Sigma = \tilde{\Sigma}.$$

It is clear that the matrix $\tilde{\Sigma}$ is positive definite. For any vector $\psi_* \in \text{int}(\Psi)$, define the sequence of random functions

$$\mathcal{T}_n(\psi_*; \mathbf{X}_n) = \partial^2 Q_n(\psi_*; \mathbf{X}_n) - \partial^2 Q_n(\psi_0; \mathbf{X}_n),$$

and denote by $\mathcal{T}_n(\psi_*; \mathbf{X}_n)_{\ell k}$ its $(\ell, k)th$ entry. Any entry of $\partial^2 Q_n(\psi_0; \mathbf{X}_n)$ is either a constant times the sum over $i = 1, \ldots, n$ of the product of the components or entries of $\partial m(\rho_0; Z_{i-1})$, $\partial^2 m(\rho_0; Z_{i-1})$, $\partial \sigma(\theta_0; Z_{i-1})$, $\partial^2 \sigma(\theta_0; Z_{i-1})$ and the random functions $\sigma(\theta_0; Z_{i-1})$, $\varepsilon_i(\psi_0)$, $\xi_i(\psi_0)$, $\dot{\xi}_i(\psi_0)$, $\zeta_i(\psi_0)$, $\dot{\zeta}_i(\psi_0)$ and $\ddot{\zeta}_i(\psi_0)$, or sums or differences of such terms. We have for example:

$$
\begin{aligned}
\partial^2 Q_n(\psi_0; \mathbf{X}_n)_{12} &= -\sum_{i=1}^n \sigma^{-1}(\theta_0; Z_{i-1})\Big(\frac{\partial^2 m(\rho_0; Z_{i-1})}{\partial \rho_1 \partial \rho_2}\xi_i(\psi_0) \\
&\quad -\sigma^{-1}(\theta_0; Z_{i-1})\frac{\partial m(\rho_0; Z_{i-1})}{\partial \rho_1}\frac{\partial m(\rho_0; Z_{i-1})}{\partial \rho_2}\dot{\xi}_i(\psi_0)\Big).
\end{aligned}
$$

In view of the assumptions $(\mathcal{A}_1)$-$(\mathcal{A}_3)$ and $(\mathcal{B}_0)$-$(\mathcal{B}_4)$, we can deduce from the above example that for each $(\ell, k)$, there exists a positive real-valued function $\mu_{\ell k}$ with $E[\mu_{\ell k}^4(Z_0)] < \infty$ such that

$$
|\mathcal{T}_n(\psi_*; \mathbf{X}_n)_{\ell k}| \leq ||\psi_* - \psi_0||_{\mathcal{E}} \sum_{i=1}^n \mu_{\ell k}(Z_{i-1}).
$$

Then, for $\delta > 0$ such that $||\psi - \psi_0||_{\mathcal{E}} < \delta$, $(n\delta)^{-1}|\mathcal{T}_n(\psi_*; \mathbf{X}_n)_{\ell k}|$ is bounded from the right by $n^{-1}\sum_{i=1}^n \mu_{\ell k}(Z_i)$ which, by the ergodic theorem, converges almost surely to $E[\mu_{\ell k}(Z_0)] < \infty$ as $n$ tends to infinity. One can thus conclude that for all $(\ell, k)$, (3.3) holds with $\mathcal{T}_n(\psi_*; \mathbf{X}_n)_{\ell k}$. Here also, as in the proof of Theorem 3.1, there exists a sequence of estimators $\hat{\psi}_n = (\hat{\rho}'_n, \hat{\theta}'_n)'$ such that, a.s., $\hat{\psi}_n \longrightarrow \psi_0$, and for any $\epsilon > 0$, there exists an event $\mathcal{S}_2$ with $P(\mathcal{S}_2) > 1 - \epsilon$, and an integer $n_2$ such that on $\mathcal{S}_2$, for $n > n_2$, $\partial Q_n(\hat{\psi}_n; \mathbf{X}_n) = \mathbf{0}$ and $Q_n(\psi; \mathbf{X}_n)$ attains a relative minimum at $\psi = \hat{\psi}_n$. Since a relative minimum for $Q_n(\psi; \mathbf{X}_n)$ is a relative maximum for $L_n(\psi; \mathbf{X}_n)$, the first part of our result is established. For the second part, it remains to prove that $n^{-1/2}\partial Q_n(\psi_0; \mathbf{X}_n)$ converges in distribution to a Gaussian random vector with mean $\mathbf{0}$ and covariance matrix $\Lambda$. This result is handled if one puts in Theorem 1 of Ngatchou-Wandji (2005) : $\omega_i = \varepsilon_i(\psi)$, $Y_i = Z_{i-1}$, $\Pi_1(z) = \sigma^{-1}(\theta_0; z)\partial m(\rho_0; z)$; $\Pi_2(z) = \sigma^{-2}(\theta_0; z)\partial \sigma(\theta_0; z)$; $\Gamma_1(x) = \phi_f(x)$; $\Gamma_2(x) = x\phi_f(x) - 1$. Finally, applying again the second part of Theorem 3.2.23 of Taniguchi and Kakizawa (2000) one has that

$$
\sqrt{n}(\hat{\psi}_n - \psi_0) \xrightarrow{\mathcal{D}} \mathcal{N}(\mathbf{0}, \Sigma^{-1}\Lambda\Sigma^{-1}).
$$

$\square$

**Corollary 4.1.** *Assume that the assumptions of Theorem 4.1 hold, and that the equalities $\int \phi_f^{(1)}(x)f(x)dx = \int \phi_f^2(x)f(x)dx$, $\int x\phi_f^2(x)f(x)dx = \int \phi_f(x)(x\phi_f(x) - 1)f(x)dx$ and $\int x(\phi_f(x) + x\phi_f^{(1)}(x))f(x)dx = \int (x\phi_f(x) - 1)^2 f(x)dx$ hold. Then*

$$
\sqrt{n}(\hat{\psi}_n - \psi_0) \xrightarrow{\mathcal{D}} \mathcal{N}(\mathbf{0}, \Sigma^{-1}), \quad n \to \infty.
$$

The conditions on the integrals in the above Theorem 4.1 and Corollary 4.1 are verified at least by Gaussian density functions $f$. When those in Corollary 4.1 are satisfied, the Fisher information matrix converges to $\Sigma$. Hence, the Cramer-Rao bound is asymptotically achieved and $\hat{\psi}_n$ is asymptotically efficient.

## 5. Some examples

Here we list some common examples that are comprised in (1.1). It is not difficult to see that the $AR(q)$ models of finite order $q$, either linear or nonlinear are within (1.1) for $\sigma(\theta; z) = Cst$. The usual ones are for example AR, SETAR, TARCH, EXPAR (see Tong (1990)). Taking $m(\rho; z) = 0$ in (1.1) yields ARCH($q$) models. For finite $q$, the most popular one is the ARCH ($q$) model obtained with

$$\sigma(\theta; Z_{i-1}) = \sqrt{\theta_0 + \theta_1 X_{i-1}^2 + \ldots + \theta_q X_{i-q}^2}, \ \ \theta_0 > 0, \ \theta_i \geq 0, \ i = 1, \ldots, q. \ \ (5.1)$$

For $q = \infty$, many other common models are within (1.1). It is the case for invertible ARMA models. In the particular case of MA(1) model defined by

$$X_i = \varepsilon_i + \theta \varepsilon_{i-1}, \ \ |\theta| < 1,$$

one has $\varepsilon_i = \sum_{j \geq 0} (-\theta)^j X_{i-j}$ from which it results that

$$X_i = \sum_{j \geq 0} (-\theta)^j X_{i-j-1} + \varepsilon_i.$$

As can be seen for instance in Peng and Yao (2003), GARCH($p, q$) models are also within (1.1). In the particular case of GARCH(1,1) model defined by

$$X_i = h_i \varepsilon_i, \ \ h_i^2 = c + a X_{i-1}^2 + b h_{i-1}^2, \ \ (5.2)$$

it is proved that for $a + b < 1$, $h_i^2 = \dfrac{c}{1-a} + b \sum_{j \geq 1} a^{j-1} X_{i-j}^2$, and consequently,

$$X_i = \left( \frac{c}{1-a} + b \sum_{j \geq 1} a^{j-1} X_{i-j}^2 \right)^{\frac{1}{2}} \varepsilon_i.$$

The class of models (1.1) for $q = \infty$ also contains invertible bilinear models, such as the subdiagonal bilinear model defined by

$$X_i = b X_{i-2} \varepsilon_{i-1} + \varepsilon_i.$$

For this model, it follows from page 103 of Taniguchi and Kakizawa (2000) that if $b^2 < 1$, then

$$\varepsilon_i = X_i + \sum_{j \geq 1} (-b)^j X_{i-j} \prod_{k=1}^{j} X_{i-k-1},$$

which in turn yields

$$X_i = bX_{i-2}\left(X_{i-1} + \sum_{j\geq 1}(-b)^j X_{i-j-1} \prod_{k=1}^{j} X_{i-k-2}\right) + \varepsilon_i.$$

**Remark 5.1.** *Although conditional least-squares and conditional maximum likelihood estimators for $\psi_0$ exists, their computation may need numerical methods, even for Gaussian density functions $f$. For example, for pure ARCH(1) models defined by (5.1) with $q = 1$ and Gaussian $f$, the conditional maximum likelihood estimators will be obtained by solving in $\psi = (\theta_0', \theta_1')'$ the equations*

$$\begin{cases} \sum_{i=1}^{n}\left(\frac{1}{\theta_0 + \theta_1 X_{i-1}^2} - \frac{X_i^2}{(\theta_0 + \theta_1 X_{i-1}^2)^2}\right) = 0 \\ \sum_{i=1}^{n}\left(\frac{X_{i-1}^2}{\theta_0 + \theta_1 X_{i-1}^2} - \frac{X_i^2 X_{i-1}^2}{(\theta_0 + \theta_1 X_{i-1}^2)^2}\right) = 0, \end{cases}$$

*with the restrictions $\theta_0 > 0$ and $0 \leq \theta_1 < 1$. This will generally need a numerical method. A similar remark can be done for the GARCH(1,1) and bilinear models defined above when estimating by either the least-squares or likelihood methods.*

## 6. Kernel estimator for the noise's density and its derivatives

In time series analysis, the conditional distribution can be very useful for the study of nonlinear phenomena such as the symmetry and the multimodality structure of a time series. In the setting of models (1.1), the conditional distribution is the distribution of the noise. Nonparametric estimation of conditional distribution has been studied among others by Hyndman and Yao (2002) who use a kernel method, Fan, Yao and Tong (1996) who use the local polynomials approach, Fan and Yim (2004) who use a cross-validation method, and Hyndman and Yao (2002) who use a kernel method and derive a test for conditional symmetry from their estimator. Bai and Ng (2001) also derive a test for conditional symmetry which rest on the kernel estimators of the conditional density and its derivatives.

In this section, we assume that the $\varepsilon_i$'s have an unknown uniformly continuous density function $f$, and we define its kernel estimator and those of its derivatives. We show the uniform consistency of these estimators. The results of this section can lead to the derivation of adaptative estimators for $\psi_0$, or to the construction of some goodness-of-fit tests for the function $f$. However, we will not study these problems here.

For all $i \in \mathbb{Z}$ and $\psi = (\rho', \theta')' \in \Psi$, we define the random function

$$\varepsilon_i(\psi) = \frac{X_i - m(\rho; Z_{i-1})}{\sigma(\theta; Z_{i-1})}. \tag{6.1}$$

Let $\psi_n = (\rho_n', \theta_n')'$ be any consistent estimator of $\psi_0$ such that $n^{1/2}(\psi_n - \psi_0)$ converges in distribution to a Gaussian distribution with mean **0** and variance matrix $\Gamma$. Take for example the least-squares estimator of Section 3, or the pseudo-maximum likelihood estimator which can be obtained from Section 4 with Gaussian $f$. Let $p$ be a nonnegative integer and $K$ be a kernel function differentiable up to order $p + 1$, with modulus of continuity $\omega_K$. Let $(h_n)$ be a sequence of real numbers such that $h_n \longrightarrow 0$, as $n$ tends to infinity. For $n = 1, 2, \ldots$, and for all $x \in \mathbb{R}$, define the random functions

$$f_n^{(p)}(\psi; x) = \frac{1}{nh_n^{p+1}} \sum_{i=1}^{n} K^{(p)} \left( \frac{x - \varepsilon_i(\psi)}{h_n} \right). \tag{6.2}$$

For observable $\varepsilon_i(\psi)$'s, the convergence of the above Bhattacharya's estimators for $f^{(p)}(x)$ is studied in Silverman (1978). Here, the $\varepsilon_i(\psi)$'s are not observable and it is natural to estimate $f^{(p)}(x)$ by $f_n^{(p)}(\psi_n; x)$. Following Singh (1979), or the more recent paper of Horová, Vieu and Zelinka (2002), other estimators of $f^{(p)}(x)$ could be defined. The results of this section are established under the following assumptions of Silverman (1978):

($\mathcal{H}_1$) $K$ is uniformly continuous with bounded variation

($\mathcal{H}_2$) $\int |K(x)|dx < \infty$ and $K(x) \longrightarrow 0$ as $|x| \to \infty$

($\mathcal{H}_3$) $\int K(x)dx = 1$

($\mathcal{H}_4$) $\int |x \log(|x|)|^{1/2} dK(x) < \infty$

($\mathcal{H}_5$) $\int_0^1 [\log(1/u)]^{1/2} d\tau(u) < \infty$, where $\tau(u) = [\omega_K(u)]^{1/2}$

($\mathcal{H}_6$) For $j = 0, \ldots, p+1$, $K^{(j)}(x) \longrightarrow 0$ as $|x| \to \infty$ and $\int |K^{(j)}(x)|dx < \infty$

($\mathcal{H}_7$) The Fourier transform of $K$ is not identically one in any neighborhood of 0.

In Silverman (1978), the assumptions ($\mathcal{H}_1$)-($\mathcal{H}_5$) are needed for the convergence of $f_n(\psi_0; x)$ to $f(x)$, while ($\mathcal{H}_1$), ($\mathcal{H}_2$), ($\mathcal{H}_4$), ($\mathcal{H}_5$)-($\mathcal{H}_7$) allow for the convergence of $f_n^{(p)}(\psi_0; x)$ to $f^{(p)}(x)$, $p \geq 1$. These assumptions hold at least for Gaussian kernels.

We have the following theorem:

**Theorem 6.1.** *Assume* ($\mathcal{B}_0$) *and* ($\mathcal{H}_6$) *hold, and the function* $K^{(p+1)}$ *is continuous. Let $r$ be any integer such that $0 \leq r \leq p$. Assume $n^{1/2} h_n^{r+2} \longrightarrow \infty$, as $n \to \infty$. Then*

$$\sup_{x \in \mathbb{R}} |f_n^{(r)}(\psi_0; x) - f_n^{(r)}(\psi_n; x)| = o_P(1).$$

*Proof.* Let $0 \leq r \leq p$, $r$ integer. By a Taylor expansion of order one, we have,

for some vector $\psi_*$ lying between $\psi_0$ and $\psi_n$:

$$
\begin{aligned}
& f_n^{(r)}(\psi_0; x) - f_n^{(r)}(\psi_n; x) \\
&= \frac{1}{nh_n^{r+1}} \sum_{i=1}^{n} \left\{ K^{(r)} \left( \frac{x - \varepsilon_i(\psi_0)}{h_n} \right) - K^{(r)} \left( \frac{x - \varepsilon_i(\psi_n)}{h_n} \right) \right\} \\
&= \frac{1}{nh_n^{r+2}} \sum_{i=1}^{n} K^{(r+1)} \left( \frac{x - \varepsilon_i(\psi_*)}{h_n} \right) \partial' \varepsilon_i(\psi_*)(\psi_0 - \psi_n) \\
&= \left\{ \frac{1}{n^{3/2}h_n^{r+2}} \sum_{i=1}^{n} K^{(r+1)} \left( \frac{x - \varepsilon_i(\psi_*)}{h_n} \right) \partial' \varepsilon_i(\psi_*) \right\} \sqrt{n}(\psi_0 - \psi_n).
\end{aligned}
$$

By the triangle inequality, it then follows that

$$
\begin{aligned}
& \sup_{x \in \mathbb{R}} |f_n^{(r)}(\psi_0; x) - f_n^{(r)}(\psi_n; x)| \\
&\quad \leq \sup_{\psi \in \Psi} \sup_{x \in \mathbb{R}} \frac{1}{nh_n^{r+2}} \sum_{i=1}^{n} \left| K^{(r+1)} \left( \frac{x - \varepsilon_i(\psi)}{h_n} \right) \right| \times \sup_{\psi \in \Psi} ||\partial' \varepsilon_i(\psi)||_{\mathcal{E}} ||\psi_0 - \psi_n||_{\mathcal{E}} \\
&\quad \leq \frac{1}{n^{3/2}h_n^{r+2}} \sum_{i=1}^{n} \left[ \sup_{\psi \in \Psi} \sup_{x \in \mathbb{R}} \left| K^{(r+1)} \left( \frac{x - \varepsilon_i(\psi)}{h_n} \right) \right| \right] \\
&\quad\quad \times \sup_{\psi \in \Psi} ||\partial' \varepsilon_i(\psi)||_{\mathcal{E}} ||\sqrt{n}(\psi_0 - \psi_n)||_{\mathcal{E}}.
\end{aligned}
$$

Since the function $K^{(r+1)}$ is continuous, by $(\mathcal{H}_6)$ it is bounded, and there exists a constant $C > 0$ such that almost surely,

$$
\sup_{\psi \in \Psi} \sup_{x \in \mathbb{R}} \left| K^{(r+1)} \left( \frac{x - \varepsilon_i(\psi)}{h_n} \right) \right| \leq C.
$$

Also, under $(\mathcal{B}_0)$ and $(\mathcal{A}_1)$, one can find a positive function $\chi(z)$ with $E[\chi^4(Z_0)] < \infty$ such that for all $1 \leq i \leq n$,

$$
\sup_{\psi = (\psi_1, \psi_2) \in \Psi} ||\partial' \varepsilon_i(\psi)||_{\mathcal{E}} \leq \chi(Z_{i-1}).
$$

From these two inequalities, we have

$$
\sup_{x \in \mathbb{R}} |f_n^{(r)}(\psi_0; x) - f_n^{(r)}(\psi_n; x)| \quad \leq \quad \frac{C}{n^{3/2}h_n^{r+2}} \left( \sum_{i=1}^{n} \chi(Z_{i-1}) \right) ||\sqrt{n}(\psi_0 - \psi_n)||_{\mathcal{E}}.
$$

By our assumptions, we have that $||\sqrt{n}(\psi_0 - \psi_n)||_{\mathcal{E}}$ converges in distribution to $||\mathcal{N}(\mathbf{0}, \Gamma)||_{\mathcal{E}}$. By the ergodic theorem, almost surely,

$$
\frac{1}{n} \sum_{i=1}^{n} \chi(Z_{i-1}) \longrightarrow E[\chi(Z_0)].
$$

The result then follows by the fact that $n^{1/2}h_n^{r+2} \longrightarrow \infty$, as $n \to \infty$. $\qquad \square$

An immediate consequence of both Theorem 4.1 and Theorems A and C of Silverman (1978), is the following corollary.

**Corollary 6.1.** *Assume that* $(\mathcal{B}_0)$ *holds and the function* $K^{(p+1)}$ *is continuous.*

(i) *Assume that* $(\mathcal{H}_1)$-$(\mathcal{H}_4)$ *hold and that* $n^{1/2}h_n^2 \longrightarrow \infty$ *and* $(nh_n)^{-1}\log(n)$ $\longrightarrow 0$ *as* $n \to \infty$. *Then uniformly, in probability,* $f_n(\psi_n; x)$ *converges to* $f(x)$.

(ii) *For* $p > 0$, *assume that the function* $K^{(p)}$ *satisfies* $(\mathcal{H}_1)$-$(\mathcal{H}_2)$, $(\mathcal{H}_4)$-$(\mathcal{H}_7)$, *and* $n^{1/2}h_n^{p+2} \longrightarrow \infty$, $n^{-1}h_n^{-2p-1}\log(1/h_n) \longrightarrow 0$ *as* $n \to \infty$. *Then uniformly, in probability,* $f_n^{(p)}(\psi_n; x)$ *converges to* $f^{(p)}(x)$.

*Proof.* By the triangle inequality, write, for $r = 0$ or $r = p$,

$$
\sup_{x \in \mathbb{R}} |f^{(r)}(x) - f_n^{(r)}(\psi_n; x)| \quad \leq \quad \sup_{x \in \mathbb{R}} |f^{(r)}(x) - f_n^{(r)}(\psi_0; x)|
$$
$$
+ \sup_{x \in \mathbb{R}} |f_n^{(r)}(\psi_0; x) - f_n^{(r)}(\psi_n; x)|,
$$

and apply Theorem 4.1 and Theorems A and C of Silverman (1978). □

**Remark 6.1.** *Take* $h_n = Cst.n^{-1/9}$. *Then, one has* $n^{1/2}h_n^2 \longrightarrow \infty$ *and* $(nh_n)^{-1}\log(n) \longrightarrow 0$ *as* $n \to \infty$, *which satisfies the hypotheses of the part* (i) *of the above corollary. For* $1 \leq p \leq 2$, *one has* $n^{1/2}h_n^{p+2} \longrightarrow \infty$ *and* $n^{-1}h_n^{-2p-1}\log(1/h_n) \longrightarrow 0$ *as* $n \to \infty$ *and the requirements of the part* (ii) *is satisfied.*

## 7. Simulation study

To illustrate some of our results, we conducted a simulation experiment that we present and comment in this last section. We restricted to models for which we could obtain explicit and simple expressions for the estimators. This avoided the use of numerical methods. We consider the following models :

$$
X_i = [\rho_0 + \rho_1 \exp(-\kappa X_{i-1}^2)]X_{i-1} + \sqrt{\theta_0 + \theta_1 X_{i-1}^2}\,\varepsilon_i, \tag{7.1}
$$

where the parameters $\rho_0, \rho_1, \kappa > 0, \theta_0 > 0$ and $\theta_1 \geq 0$ eventually satisfy some conditions insuring the existence, the invertibility, the stationarity and the ergodicity of $(X_i)_{i \in \mathbb{Z}}$. For example, for model (ii) below, (7.1) admits a strictly stationary and geometrically ergodic solution $(X_i)_{i \in \mathbb{Z}}$ as soon as $0 \leq \theta_1 < 1$. The noise densities $f$ that we used were either Gaussian or Laplace. More precisely, we studied the cases

(i) $\rho_0 = 0$, $0 < \rho_1 < 1$, $\kappa = 0.1$ and $\theta_1 = 0$, with $f$ either Gaussian or Laplace.
(ii) $\rho_0 = 0$, $\rho_1 = 0$, $\kappa = 0$ and $0 < \theta_1 < 1$, with $f$ Gaussian.
(iii) $\rho_1 = 0$, $\kappa = 0$ and $0 < \theta_1 < 1$, with $f$ Gaussian.

Except model (i) with Gaussian $f$ and model (ii), there is no guaranty that $(X_i)_{i \in \mathbb{Z}}$ be stationary and / or ergodic for the other models.

For the computation of least-squares estimators, the weight functions were $\lambda(z) = \gamma(z) \equiv 1$, which yield the classical least-squares estimators. In each case, our estimates were computed on the basis of 1,000 samples of length $n$. For model (i), from simple computations, it is easy to see that the least-squares estimator coincides with the maximum likelihood estimator for Gaussian $f$. The results concerning this model are listed in Table 1. They show that, for samples of size $n = 50$, and for either density considered, the least-squares estimators are close to the true value of the parameters. Rapid calculus show that these estimators are unbiased. The trials for this model were also done for $n \geq 100$. From the results that we do not present here, the estimates obtained were more accurate. Concerning the models (ii), the results were in general better for the maximum likelihood estimators than least-squares, for all the sample sizes $n = 100, 200, 400$ (see Tables 2 and 3). Both estimators moved to the true parameter as $n$ grew. For the models (iii), only least-squares estimators were computed. This was done for $n = 100, 200, 400$. We observed in these cases that the estimates of $\rho_0$ were good while $\theta_0$ was always overestimated and $\theta_1$ was underestimated (see Table 4). The least-squares estimates for the models (ii) also behaved this way. This likely comes from the fact that the least-squares estimators for these models are highly biased. It seems from our simulation experiment that their bias converge slowly to 0, as $n$ grows.

TABLE 1
*Conditional least-squares estimator for the parameters of model (i) with Gaussian noises (middle two columns) and Laplace noises (last two columns) and sample size n = 50*

| $\rho_1$ | $\theta_0$ | | $\tilde{\rho}_1$ | $\tilde{\theta}_0$ | | $\tilde{\rho}_1$ | $\tilde{\theta}_0$ |
|---|---|---|---|---|---|---|---|
| -0.80 | 0.10 | | -0.77 | 0.098 | | -0.778 | 0.097 |
| -0.50 | 0.50 | | -0.47 | 0.493 | | -0.482 | 0.486 |
| 0.20 | 0.10 | | 0.191 | 0.098 | | 0.197 | 0.098 |
| 0.80 | 0.80 | | 0.763 | 0.795 | | 0.782 | 0.779 |
| 0.90 | 1.00 | | 0.864 | 0.998 | | 0.886 | 0.976 |

TABLE 2
*Conditional least-squares estimator for the parameters of ARCH(1) model (ii) ($\rho_0 = \rho_1 = 0$) and Gaussian noises, for sample sizes n = 100, n = 200 and n = 400*

| n= | | | 100 | | 200 | | | 400 | |
|---|---|---|---|---|---|---|---|---|---|
| $\theta_0$ | $\theta_1$ | | $\tilde{\theta}_0$ | $\tilde{\theta}_1$ | $\tilde{\theta}_0$ | $\tilde{\theta}_1$ | | $\tilde{\theta}_0$ | $\tilde{\theta}_1$ |
| 0.40 | 0.30 | | 0.433 | 0.210 | 0.428 | 0.243 | | 0.418 | 0.264 |
| 0.50 | 0.20 | | 0.525 | 0.146 | 0.515 | 0.169 | | 0.510 | 0.180 |
| 0.30 | 0.10 | | 0.304 | 0.075 | 0.304 | 0.089 | | 0.300 | 0.095 |
| 0.40 | 0.40 | | 0.472 | 0.271 | 0.451 | 0.306 | | 0.442 | 0.329 |
| 0.60 | 0.05 | | 0.610 | 0.031 | 0.608 | 0.036 | | 0.603 | 0.044 |

TABLE 3
*Conditional maximum likelihood estimator for the parameters of ARCH(1) model (ii) ($\rho_0 = \rho_1 = 0$) and Gaussian noises, for sample sizes n = 100, n = 200 and n = 400*

| n= | | | 100 | | 200 | | | 400 | |
|---|---|---|---|---|---|---|---|---|---|
| $\theta_0$ | $\theta_1$ | | $\hat{\theta}_0$ | $\hat{\theta}_1$ | $\hat{\theta}_0$ | $\hat{\theta}_1$ | | $\hat{\theta}_0$ | $\hat{\theta}_1$ |
| 0.40 | 0.30 | | 0.413 | 0.268 | 0.407 | 0.284 | | 0.401 | 0.297 |
| 0.50 | 0.20 | | 0.508 | 0.175 | 0.505 | 0.188 | | 0.503 | 0.191 |
| 0.30 | 0.10 | | 0.294 | 0.105 | 0.300 | 0.100 | | 0.299 | 0.100 |
| 0.40 | 0.40 | | 0.415 | 0.364 | 0.406 | 0.381 | | 0.402 | 0.393 |
| 0.60 | 0.05 | | 0.583 | 0.067 | 0.593 | 0.055 | | 0.596 | 0.051 |

TABLE 4
*Conditional least-squares estimator for the parameters of model (iii) with $\rho_1 = 0$ and
Gaussian noises, for sample sizes $n = 100$, $n = 200$ and $n = 400$*

| n= | | | | 100 | | | | 200 | | |
|---|---|---|---|---|---|---|---|---|---|---|
| $\rho_0$ | $\theta_0$ | $\theta_1$ | | $\tilde{\rho}_0$ | $\tilde{\theta}_0$ | $\tilde{\theta}_1$ | | $\tilde{\rho}_0$ | $\tilde{\theta}_0$ | $\tilde{\theta}_1$ |
| 0.20 | 0.40 | 0.30 | | 0.189 | 0.447 | 0.184 | | 0.196 | 0.435 | 0.219 |
| 0.30 | 0.50 | 0.20 | | 0.292 | 0.534 | 0.131 | | 0.292 | 0.524 | 0.155 |
| 0.50 | 0.30 | 0.10 | | 0.491 | 0.313 | 0.060 | | 0.494 | 0.307 | 0.075 |
| 0.60 | 0.40 | 0.05 | | 0.582 | 0.411 | 0.025 | | 0.591 | 0.408 | 0.033 |
| 0.40 | 0.40 | 0.10 | | 0.390 | 0.415 | 0.062 | | 0.389 | 0.410 | 0.077 |

| n= | | | | 400 | | |
|---|---|---|---|---|---|---|
| $\rho_0$ | $\theta_0$ | $\theta_1$ | | $\tilde{\rho}_0$ | $\tilde{\theta}_0$ | $\tilde{\theta}_1$ |
| 0.20 | 0.40 | 0.30 | | 0.198 | 0.424 | 0.253 |
| 0.30 | 0.50 | 0.20 | | 0.295 | 0.517 | 0.173 |
| 0.50 | 0.30 | 0.10 | | 0.494 | 0.303 | 0.085 |
| 0.60 | 0.40 | 0.05 | | 0.594 | 0.405 | 0.040 |
| 0.40 | 0.40 | 0.10 | | 0.399 | 0.405 | 0.086 |

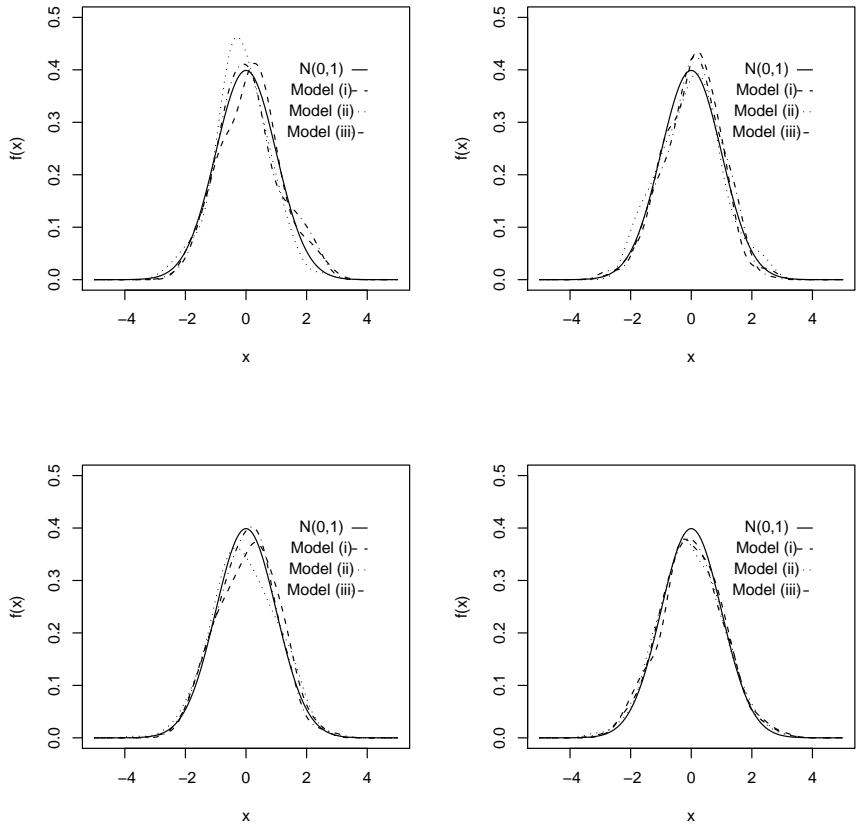

FIG 1. *Estimation of the density of the noise in models (i), (ii) and (iii) for $n = 100, 200, 400, 600$.*

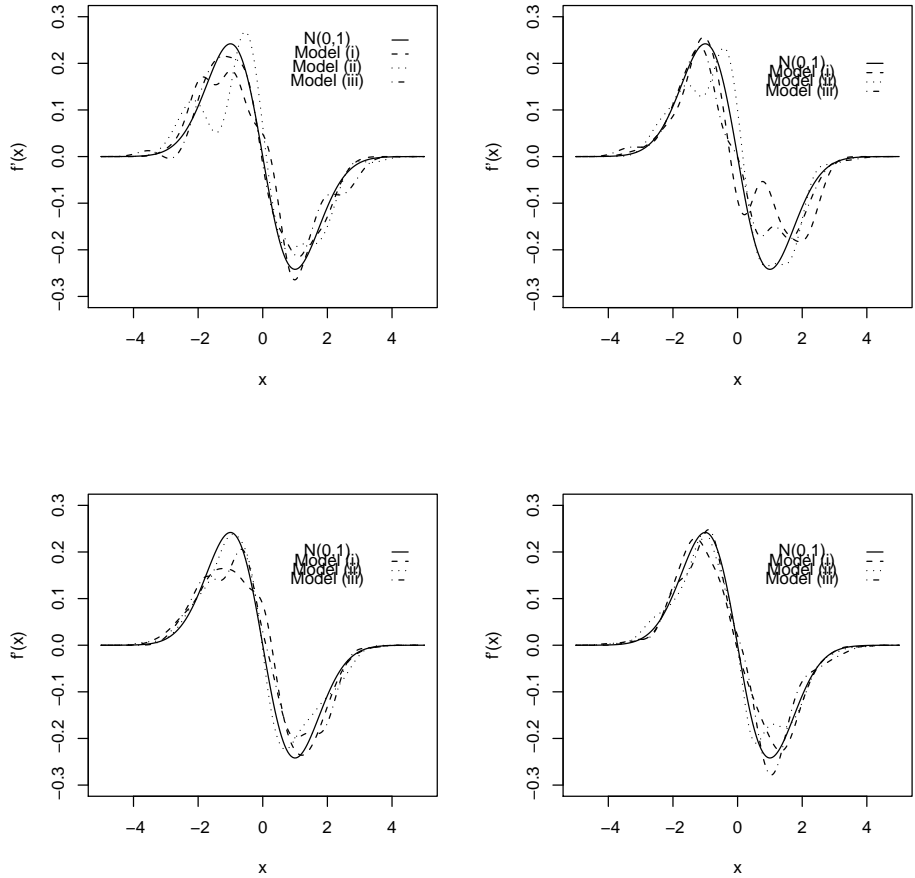

FIG 2. *Estimation of the first derivative of the density of the noise in models (i), (ii) and (iii) for* $n = 100, 200, 400, 600$

For the kernel density estimation, we restricted our trials to models (*i*)-(*iii*) with Gaussian density. The residuals were computed from least-squares fit with $\gamma(z) = \lambda(z) \equiv 1$. We took $\psi_n = \tilde{\psi}_n$ (see Sections 3 and 5), a Gaussian kernel with $h_n = c_n n^{-1/9}$, where, denoting by $\sigma_n$ the empirical standard deviation and $X_{n,\frac{1}{4}}$ and $X_{n,\frac{3}{4}}$ the first and third empirical quartiles of $(X_1, \ldots, X_n)$,

$$c_n = \frac{0.9 \min\{\sigma_n, (X_{n,\frac{3}{4}} - X_{n,\frac{1}{4}})\}}{1.34}.$$

This sequence $(c_n)$ given in the software R seemed to give better results than $c_n = \sigma_n$. It is easy to check that the Gaussian kernel and the sequences $(h_n)$ clearly satisfy the assumptions of Theorem 6.1. We took $\rho_1 = -0.5$, $\theta_1 = 1$ for

($i$), $\theta_0 = 0.4$, $\theta_1 = 0.1$ for ($ii$) and $\rho_0 = 0.6$, $\theta_0 = 0.4$, $\theta_1 = 0.05$ for ($iii$). The different plots of $\tilde{f}_n$ and $f$ are gathered on the same graph. The same is done for $\tilde{f}_n^{(1)}$ and $f^{(1)}$. The trials were done for $n = 100, 200, 400, 600$. The estimates obtained for the density were good (see Figure 1). Those of the derivative of the density were not good for $n = 100, 200$, especially in the vicinity of the maxima. They were better for $n = 400, 600$ (see Figure 2). For the density and its first derivative, one can see that the estimates from the models ($ii$) and ($iii$) were not very close to the true functions. This is probably due to the sampling fluctuations or to the bias of the least-squares estimators of the parameters of these models. The good behavior of the estimates obtained from ($i$) may come from the fact that the conditional likelihood and the conditional least-squares estimators of the parameter $\psi_0$ are the same in this case as we earlier pointed out.